\documentclass[12pt]{article}
\usepackage{amsmath}
\usepackage{amssymb}

\newcommand{\R}{{\mathbb R}}

\newcommand{\fhat}{\hat{f}}
\newcommand{\fn}{\!:\!}
\newcommand{\lsum}{\sum\limits}

\newcommand{\lint}{\int\limits}

\newcommand{\rll}{Riemann--Lebesgue Lemma}
\newcommand{\qed}{\mbox{$\quad\blacksquare$}}

\begin{document}
\begin{center}
{\large\bf Rapidly growing Fourier integrals}
\vskip.25in
Erik Talvila\\ [2mm]
\end{center}
{\bf 1. The Riemann--Lebesgue Lemma.}\quad 
In its usual form, the \rll\ reads as follows:  If $f\in L^1$ and 
$\hat{f}(s)=\int_{-\infty}^{\infty}e^{isx}f(x)\,dx$ is its Fourier transform,
then $\hat{f}(s)$ exists and is finite for each  $s\in\R$ and 
$\hat{f}(s)\to 0$ as $|s|\to\infty$ ($s\in\R$).  
This result encompasses Fourier sine
and cosine transforms  as well as
Fourier series coefficients for functions periodic on finite intervals.  
When the integral
is allowed to converge conditionally, the asserted asymptotic behaviour can
fail dramatically.
In fact, we show that for each sequence $a_n\uparrow\infty$ we can
find a continuous function $f$ such that $\fhat(s)$ exists for each $s\in\R$ 
and
$\fhat(n)\geq a_n$ for all integers $n\geq 1$.  We also work out
the asymptotics of a class of Fourier integrals that can have arbitrarily
large polynomial growth.  Our main tool is the principle of stationary
phase.  The conditionally convergent integrals we consider in this paper can be 
thought of as Henstock integrals \cite{bartle} or as improper Riemann
integrals.

Two examples of conditionally convergent 
Fourier transforms that do not tend to zero at infinity
can be obtained from \cite[3.691]{gradshteyn}:
\begin{equation}
\lint_{x=0}^{\infty}\left\{\!\!\begin{array}{c}
\sin(ax^2)\\
\cos(ax^2)
\end{array}
\!\!\right\}
\cos(sx)\,dx  = 
\frac{1}{2}\sqrt{\frac{\pi}{2a}}\left[\cos\left(\frac{s^2}{4a}\right)
\mp\sin\left(\frac{s^2}{4a}\right)\right]\label{I.4}
\end{equation}
and
\begin{eqnarray}
\lefteqn{\lint_{x=0}^{\infty}\left\{\!\!\begin{array}{c}
\sin(ax^2)\\
\cos(ax^2)
\end{array}
\!\!\right\}
\sin(sx)\,dx  =} \notag\\
 & & \sqrt{\frac{\pi}{2a}}
 \left[\left\{\!\!\begin{array}{c}\cos\left[s^2/(4a)\right]\\
 \sin\left[s^2/(4a)\right]
 \end{array}
 \!\!\right\}
C\left(\frac{s^2}{4a}\right)\pm
 \left\{\!\!\begin{array}{c}\sin\left[s^2/(4a)\right]\\
 \cos\left[s^2/(4a)\right]
 \!\!\end{array}
 \right\}
S\left(\frac{s^2}{4a}\right)\right].
 \label{I.5}
 \end{eqnarray}
Here, $a>0$ and 
\begin{equation}
C(x)=  
\frac{1}{\sqrt{2\pi}}\int_{0}^x\!\!\cos t\frac{dt}{\sqrt{t}}\quad\mbox{ and }
\quad 
S(x)=\frac{1}{\sqrt{2\pi}}\int_{0}^x\!\!\sin t\frac{dt}{\sqrt{t}}
\end{equation}
are the Fresnel integrals.  Using \eqref{I.4} and \eqref{I.5} we have
the Fourier transforms of $x\mapsto\sin(ax^2)$ and $x\mapsto\cos(ax^2)$.
Both of these transforms oscillate rapidly at infinity with amplitude
that is asymptotically constant.  Note that
$C(\infty)+iS(\infty)=e^{i\pi/4}/\sqrt{2}$.  The values of all of these
integrals
are consequences of the formula $\int_{-\infty}^{\infty}e^{-x^2}\,dx
=\sqrt{\pi}$.

On a finite interval $[a,b]$ the \rll\ for conditionally convergent
integrals takes the following form:  Suppose 
$\int_a^b f$ exists.  Then
its Fourier transform, $\hat{f}(s)=\int_a^b e^{isx}f(x)\,dx$, exists 
for all $s\in
\R$ since for each $s\in\R$, the
exponential function in the integrand is
of bounded variation on the finite interval $[a,b]$.
Let $F(x)=\int_a^x f(t)\,dt$
and integrate by parts:
$$
\fhat(s)=\int_a^b e^{isx}f(x)\,dx=e^{isb}F(b)-is\int_a^b e^{isx}F(x)\,dx.
$$
Since $F\in L^1$ the \rll\ gives $\hat{f}(s)=o(s)$ as $|s|\to\infty$ in $\R$.
In \cite{tit}, Titchmarsh proved this
was the best possible estimate.  

\bigskip
\noindent
{\bf 2. Arbitrarily large pointwise growth.}\quad On the real line we 
have the following example of arbitrarily large pointwise growth.\\

\noindent
{\bf Proposition:}\quad
{\em Given any sequence of positive real numbers $\{a_n\}$,
there is a continuous
function $f\fn \R\to\mathbb{C}$ such that $\fhat(s)$ exists
for each $s\in\R$ and $\fhat(n)\geq a_n$ for all $n\geq 1$.}\\

\bigskip
{\bf Proof:} Let $\alpha_n=2a_n+1$.  We can assume that $a_n\geq 1$. 
Let $f_n(x)=\alpha_ne^{-inx}\sin(r_n x)\chi_{[-b_n,b_n]}(x)/x$ 
for $x\not=0$ and $f_n(0)=\alpha_nr_n$,
where $0<r_n\leq 1/2$ are chosen such that $\sum \alpha_nr_n<\infty$.
The sequence $b_n>0$ is to be determined so that $b_nr_n$ is an
integer multiple of $\pi$.  
Let $f(x)=(1/\pi)\sum f_n(x)$.

To compute the Fourier
transform of $f$, use the formula
\begin{equation}
\int_0^\infty \sin(ax)\frac{dx}{x}=\left\{\begin{array}{cl}
\pi/2, & a>0\\
0, & a=0\\
-\pi/2, & a<0.
\end{array}
\right.\label{II.1}
\end{equation}
This integral can be evaluated with contour integration 
\cite[p.~184]{spiegel}, with uniform convergence 
\cite[p.~262]{rogers},
and with the \rll\ \cite[p.~589]{courant}.

The estimate $|f_n(x)|\leq \alpha_ nr_n$ shows that
$\sum\alpha_nf_n(x)$ converges uniformly on $\R$.
We can
interchange orders of summation
and integration in the calculation of $\hat f$.  Let $m\geq 1$.  We then have
\begin{eqnarray*}
\fhat(m) & = & \frac{1}{\pi}\lsum_{n=1}^\infty\alpha_n\int_{-b_n}^{b_n}
e^{i(m-n)x}\sin(r_n x)\,\frac{dx}{x}\\
 & = & \frac{1}{\pi}\lsum_{n=1}^\infty \alpha_n\!\int_{0}^{b_n}\!\!
\left\{\sin\left[(|m-n|+r_n)x\right]-
\sin\left[(|m-n|-r_n)x\right]\right\}
\frac{dx}{x}.
\end{eqnarray*}
Integration by parts shows that
\begin{equation}
\left|\int_x^\infty e^{iax'}\tfrac{dx'}{x'}\right|\leq \frac{2}{ax}
\quad\text{for 
all }a, x>0.\label{II.2}
\end{equation}
Therefore, \eqref{II.1} and \eqref{II.2} ensure that
\begin{eqnarray}
\fhat(m) & \geq &\alpha_m-\frac{4\,\alpha_m}{\pi r_mb_m}
-\frac{2}{\pi}\lsum_{n\not=m}\frac{\alpha_n}{b_n}\left[\frac{1}{|m-n|+r_n}
+\frac{1}{|m-n|-r_n}\right]\notag\\
 & \geq & \alpha_m -\frac{4\,\alpha_m}{\pi r_mb_m}
 -\frac{6}{\pi}\lsum_{n=1}^\infty\frac{\alpha_n}{b_n}\notag\\
 & \geq & a_m\mbox{ if } r_mb_m\geq 12/\pi \mbox{ and } 
 b_n\geq 6(2a_n+1)2^n/\pi \mbox{ for all } n\geq 1.\label{II.2.5}
\end{eqnarray}
If we take $b_n=\lceil 6(2a_n+1)2^n/\pi\rceil\pi/r_n$ then
the conditions 
in \eqref{II.2.5}
are satisfied and $f$ is continuous on $\R$.

The interchange of summation and integration can be justified as follows:
Let $s\in\R$ and $x\geq 1$.  Define
$$
G(x)=\lsum_{n\geq |s|+1}\int_{0}^x e^{isx'}f_n(x')\,dx'.
$$
Since we already know that 
$\sum\alpha_nf_n(x)$ converges uniformly and that
$|\int_0^\infty f_n|= |\int_0^{b_n}f_n|\leq \alpha_nr_nb_n$, 
our interchange of
summation and integration is valid provided that
$G$ converges uniformly on
$[0,\infty]$ \cite[Exercise 5 in \S5.6]{rogers}.  
The Frullani integral formula \cite[p. 263]{rogers} says that
\begin{equation}
\int_0^{\infty}\left[\cos(px)-\cos(qx)\right]\frac{dx}{x}=\log\left(
\frac{q}{p}\right);\quad p,q>0.\label{II.3}
\end{equation}
We have $G(x) =[S_1(x)+S_2(x)]/2$ where
\begin{eqnarray*}
S_1(x)& = & \lsum_{\stackrel{b_n< x}{n\geq |s|+1}}\alpha_n\int_{0}^{b_n}
\left\{\sin\left[(|s-n|+r_n)x'\right]-
\sin\left[(|s-n|-r_n)x'\right]\right\}\\
 & &
\quad\quad\quad -i
\left\{\cos\left[(|s-n|+r_n)x'\right]-
\cos\left[(|s-n|-r_n)x'\right]\right\}
\,\frac{dx'}{x'}
\end{eqnarray*}
and
\begin{eqnarray*}
S_2(x)& = & \lsum_{\stackrel{b_n\geq x}{n\geq |s|+1}}\alpha_n\int_{0}^{x}
\left\{\sin\left[(|s-n|+r_n)x'\right]-
\sin\left[(|s-n|-r_n)x'\right]\right\}\\
 & &
 \quad\quad\quad -i
 \left\{\cos\left[(|s-n|+r_n)x'\right]-
 \cos\left[(|s-n|-r_n)x'\right]\right\}
 \,\frac{dx'}{x'}.
 \end{eqnarray*}
Using \eqref{II.1}, \eqref{II.2}, \eqref{II.3}, and the inequality
$0\leq\log(1+t)\leq t$ ($t\geq 0$) we obtain
\begin{eqnarray*}
S_1(x)& = & -\lsum_{\stackrel{b_n< x}{n\geq |s|+1}}\alpha_n\int_{b_n}^{\infty}
\left\{\sin\left[(|s-n|+r_n)x'\right]-
\sin\left[(|s-n|-r_n)x'\right]\right\}\\
 & &
 \quad\quad\quad -i
 \left\{\cos\left[(|s-n|+r_n)x'\right]-
 \cos\left[(|s-n|-r_n)x'\right]\right\}
 \,\frac{dx'}{x'}\\
  & & \quad\quad\quad-i\,\alpha_n\log\left[\frac{|s-n|+r_n}{|s-n|-r_n}\right]
 \end{eqnarray*}
and
\begin{eqnarray*}
|S_1(x)| & \leq & \lsum_{n\geq |s|+1}\frac{4\alpha_n}{b_n}
\left[\frac{1}{|s-n|+r_n}
+\frac{1}{|s-n|-r_n}\right] + \frac{2\alpha_nr_n}{|s-n|-r_n}\\
 & \leq & 4\lsum_{n=1}^\infty \alpha_n\left(3/b_n+r_n\right)\\
 & < & \infty.
\end{eqnarray*}
Also,
$$
S_2(x) = i \lsum_{\stackrel{b_n\geq x}{n\geq |s|+1}}\alpha_n
\left\{I(x)+\log\left[(|s-n|+r_n)/(|s-n|-r_n)\right]\right\}
$$
where
$$
I(x)=\int_{x'=x}^\infty 
\left[e^{i(|s-n|+r_n)x'}-e^{i(|s-n|-r_n)x'}\right]\frac{dx'}{x'}.
$$
Integrating by parts twice gives
\begin{eqnarray*}
I(x) & = & -\frac{1}{ix}\left[\frac{e^{i(|s-n|+r_n)x}}{|s-n|+r_n}
-\frac{e^{i(|s-n|-r_n)x}}{|s-n|-r_n}\right]\\
 & & \quad - \frac{1}{x^2}\left[\frac{e^{i(|s-n|+r_n)x}}{(|s-n|+r_n)^2}
 -\frac{e^{i(|s-n|-r_n)x}}{(|s-n|-r_n)^2}\right]\\
  & & \quad -2\int_{x'=x}^\infty\left[\frac{e^{i(|s-n|+r_n)x'}}{(|s-n|+r_n)^2}
  -\frac{e^{i(|s-n|-r_n)x'}}{(|s-n|-r_n)^2}\right]\frac{dx'}{(x')^3}.
\end{eqnarray*}
The Mean Value Theorem now shows that there are constants $c_1$, $c_2$, 
and $c_3$,
independent of $n\geq |s|+1$ and $x\geq 1$, such that
$$
|I(x)|  \leq  c_1r_n+c_2r_n\!\!\lint_{x'=x}^\infty \frac{dx'}{(x')^2}
\leq  c_3r_n.
$$
It now follows that $|S_2(x)|\leq \sum\alpha_n(c_3r_n+4r_n)<\infty$.
Hence, $G$ converges uniformly on $[0,\infty]$.
As there is a similar calculation for $x<-1$, our commutation of
$\int$ and $\sum$ in the calculation of $\fhat$
is valid.  This also shows that $\fhat(s)$ exists for every $s\in\R$.
\qed

The example can be modified so that $f$ is real-valued if we use
$f_n(x)=\alpha_n\cos(nx)\sin(r_n x)\chi_{[-b_n,b_n]}(x)/x$.  
With essentially 
the same proof we can have $\fhat(\nu_n)\geq a_n$ for
any sequence with $\nu_n-\nu_{n+1}\geq \delta$ for all
$n\geq 1$ and some $\delta>0$.
And,  if instead of the characteristic function $\chi_{[-b_n,b_n]}$ we put in a
$C^\infty$ cutoff function and take $r_n$ small enough, then $f$ can be
taken to be 
$C^\infty$ with $\fhat(n)\geq a_n$ for all $n\geq 1$.
This is very different from the Lebesgue case.  When $f\in L^1$, the smoother
$f$ is the
more rapidly $\hat{f}$ decays;  see
\cite[\S 3.4]{hardy}, and \cite[p.~45]{zygmund}.
By contrast, with conditional convergence, even for smooth $f$ we can have $\hat{f}$
growing  at an arbitrarily large rate.  Convergence of
$\int_{-\infty}^\infty f$ is necessary but not sufficient for the 
existence of $\hat{f}$.

The usual heuristic explanation of the $L^1$ \rll\ is that the
rapid oscillation of $e^{isx}$ for large
$|s|$ makes the positive and negative parts nearly cancel out in the integral
for $\hat{f}$.  Following the argument in \cite[p.~98]{krantz},
shows that $L^1$ functions are well approximated by continuous 
functions that are themselves nearly 
constant on small intervals.
For small $\epsilon>0$, there is a continuous function $g$ such that
\begin{eqnarray*}
\int_{a-\epsilon}^{a+\epsilon}e^{isx}f(x)\,dx  & \approx  &
\int_{a-\epsilon}^{a+\epsilon}e^{isx}g(x)\,dx
  \approx  g(a)\int_{a-\epsilon}^{a+\epsilon}e^{isx}\,dx\\
 & = &  2\,g(a)e^{isa}\sin(s\epsilon)/s
  \to  0\quad\text{as } |s|\to\infty.
\end{eqnarray*}
Summing such terms then shows that $\hat{f}\to 0$ as $|s|\to\infty$.
However, with conditionally convergent integrals the integrand can oscillate
with nearly the same period as $e^{isx}$ over large intervals.  For example,
in \eqref{I.4} put $a=1$.  The integrand is then $\cos(x^2)\cos(sx)$.  When
$x$ is close to $s$ the integrand is close to $\cos^2(x^2)$, which no longer
oscillates.  Thus, integrating near $s$ contributes a relatively large amount
to the integral so that $\hat{f}$ does not go to $0$ as $s\to\infty$.  Examples
in the next section also illustrate this point.

\bigskip
\noindent
{\bf 3. More Fourier integrals.}\quad
Integrating by parts shows that the integral
\begin{equation}
J=\int_{x=0}^\infty x^{\alpha}e^{i[ax^2-sx]}\,dx\label{IV.1}
\end{equation}
converges for $|\alpha|<1$.  Changing the sign of $s$ gives a similar
integral.  (Letting $\alpha\to 1^-$ then leads to two divergent integrals.
In 
\cite[3.851]{gradshteyn}
and
\cite[2.5.22]{prudnikov} they are listed as converging!  
See \cite{talvila} for a discussion of these
divergent integrals.)
Our main goal
here is to see how $J$ behaves as $s\to\infty$.  The integrand can grow nearly linearly.  How does this affect the growth of $J$?

Assume that $0<\alpha<1$.  Use the transformation $x\mapsto sx/a$.  Then
$$
J=\left(\frac{s}{a}\right)^{\alpha+1}\int_{x=0}^\infty x^{\alpha}
e^{it[(x-1/2)^2-1/4]}\,dx,
$$
where $t=s^2/a$.  As $t\to\infty$ the exponential term oscillates rapidly
except near the minimum of $x^2-x$.  Thus, we expect nearly perfect cancellation
except near $x=1/2$, and we expect that the 
major contribution to $J$ should come from integrating
near $1/2$.  Hence,
\begin{eqnarray}
J & \sim & 
\left(\frac{s}{a}\right)^{\alpha+1}2^{-\alpha}e^{-it/4}\int_{x=1/2-\epsilon}^{1/2+\epsilon}
e^{it(x-1/2)^2}\,dx\quad (t\to\infty)\notag\\
 & \sim & \left(\frac{s}{a}\right)^{\alpha+1}2^{-\alpha}e^{-it/4}\int_{x=-\infty}^{\infty}
e^{itx^2}\,dx\quad (t\to\infty)\notag\\
 & = & \left(\frac{s}{a}\right)^{\alpha+1}2^{-\alpha}e^{-it/4}\sqrt{\pi/t}\,e^{i\pi/4}
\quad (s^2/a\to\infty)\notag\\
 & = & \sqrt{\pi}\,2^{-\alpha}e^{i(\pi-s^2/a)/4}a^{-(\alpha+1/2)}\,s^{\alpha}
\quad (s^2/a\to\infty).\label{IV.2}
\end{eqnarray}
Integrating by parts on the complement of $(1/2-\epsilon,1/2+\epsilon)$ and
using the \rll\ shows that \eqref{IV.2} gives the dominant behaviour of $J$.
This heuristic argument is made precise in \cite[pp.~76--84]{wong}
and is known as the {\it principle of stationary phase}.

Fixing $a>0$ and taking $\alpha$ close to $1$ shows that the best estimate
is $J=O(s^w)$ $(s\to\infty)$, where $w<1$ but $w$ can be made arbitrarily
close to $1$.  The asymptotic behaviour of $J$ as $\alpha\to 1^-$ and
$s\to\infty$ is much more complicated.  Such uniform asymptotic approximations
are discussed in Chapter~VII of  \cite{wong}.

The integral $J$ can be evaluated in terms of
confluent hypergeometric functions \cite[pp.~23, 24, 136]{ober1},
parabolic cylinder functions \cite[pp.~23, 136]{ober2} and hypergeometric 
functions 
\cite[p.~430]{prudnikov}.
Changing the sign of $s$ puts the minimum of the exponent $ax^2+sx$ outside
the integration interval.  Integration by parts then shows that
$\int_0^\infty x^{\alpha}e^{i[ax^2+sx]}\,dx=o(s^{\alpha -1})$ ($s\to\infty$).
Linear combinations of this integral and \eqref{IV.1} lead to four integrals
akin to \eqref{I.4} and \eqref{I.5}, whose behaviour as $s\to\infty$ is given
by \eqref{IV.2}.

Let's consider one final integral.  
Let
\begin{eqnarray}
K & = & \int_{x=0}^\infty x^{\alpha}e^{i[ax^\nu-sx]}\,dx\notag\\
 & = & \left(\frac{s}{a}\right)^{\tfrac{\alpha+1}{\nu-1}}
\int_{x=0}^\infty x^{\alpha}e^{it[x^\nu-x]}\,dx \quad 
\left(t=(s^{\nu}/a)^{1/(\nu-1)}\right).\label{IV.3}
\end{eqnarray}
Integration by parts shows that \eqref{IV.3} converges for $-1<\alpha<\nu -1$.  Write
$\phi(x)=x^\nu-x$.  If $\nu>1$ then $\phi$ has a minimum at $\nu^{-1/(\nu-1)}$.
Expanding near this point gives
$$
K\sim \sqrt{\tfrac{2\pi}{\nu-1}}\,e^{i\pi/4}e^{i(s^\nu/a)^{1/(\nu-1)}
[\nu^{-\nu/(\nu-1)}-\nu^{-1/(\nu-1)}]}
\nu^{-\frac{2\alpha+1}{2(\nu-1)}}\,
a^{-\frac{2\alpha+1}{2(\nu-1)}}\,s^{\frac{2\alpha+2-\nu}{2(\nu-1)}}
$$
as $t\to\infty$.  What values of $\alpha$ and $\nu$ make $K$ large when $s\to\infty$?
If $\alpha>\nu/2 -1$, the growth of $K$ is largely 
determined by how close $\nu$ is to $1^+$.
Fix $\alpha\in(-1/2,0]$ and fix $a>0$.  
It is then apparent that the best estimate of $K$ 
as $s\to\infty$ is $O(s^w)$, where $w=(\alpha+1-\nu/2)/(\nu-1)$.
This exponent can be made arbitrarily large by taking $\nu$ close to $1^+$.
Hence, large growth in $K$ does not come from taking $\alpha$ close to
$\nu-1$ to  make the 
term $x^{\alpha}$ as  large as possible as $x\to\infty$;
rather, it
comes from flattening out the minimum of $\phi$ by making $\phi$ nearly linear.

{\it
University of Alberta,
Edmonton AB Canada T6G 2E2\\
etalvila@math.ualberta.ca}
\end{document}